\documentclass[12pt]{article}
\usepackage{amsmath, amssymb, amsthm, amsfonts, enumerate}
\usepackage{graphicx}
\usepackage{tikz}
\usepackage{caption}
\usepackage{float}
\usetikzlibrary{arrows.meta}

\newtheorem{theorem}{Theorem}[section]
\newtheorem{lemma}[theorem]{Lemma}
\newtheorem{proposition}[theorem]{Proposition}
\newtheorem{remark}[theorem]{Remark}
\newtheorem{remarks}[theorem]{Remarks}
\newtheorem{cor}[theorem]{Corollary}

\def\rr{{\mathbb R}}
\def\sik{{\rr}^2}
\def\qq{{\mathbb Q}}
\def\ff{{\cal F}}
\def\gg{{\cal G}}
\def\hh{{\cal H}}
\def\pp{{\cal P}}
\def\ttt{{\cal T}}
\def\al{\alpha}
\def\be{\beta}
\def\ga{\gamma}
\def\de{\delta}
\def\la{\lambda}

\def\ep{\varepsilon}

\def\cd{\cdot}
\def\stb{,\ldots ,}
\def\su{\subset}
\def\se{\setminus}

\def\proof{\noindent {\bf Proof.} }
\def\sumin{\sum_{i=1}^n}
\def\sumik{\sum_{i=1}^k}
\def\nl{[0,1]}

\begin{document}
\title{Invariants and equidecomposability in rings of polygons with sides
of given directions}

\author{G. Kiss and M. Laczkovich}

\maketitle

\begin{abstract}
We investigate equidecomposability in the ring of polygons with sides restricted to given directions and using only translations. Extending classical results of Dehn and Hadwiger, we prove that  equidecomposability in these rings is equivalent to the equality of some translation invariant simple valuations. We also consider the algebraic structure of direction sets.
We show that under mild conditions, equidecomposability with respect to a set $S$ of slopes of the given directions is equivalent to equidecomposability with respect to the field generated by $S$.     
We also provide a complete description of all invariants of these polygon rings.
\end{abstract}


\vspace{0.5em}
\noindent\textbf{MSC Classification (2020):} Primary: 52B45, Secondary: 51M20

\section{Introduction and main results}
By a classical theorem of M. Dehn, two rectangles with sides
parallel to the axes are equidecomposable with rectangular pieces and using
translations if and only if the area of the rectangles are equal, and if
the ratio of the vertical sides of the rectangles is rational. (See \cite{D}
and \cite[Korollar I, p.~77]{H}.)

In this note we consider the following, more general situation. By a
{\it polygon} we mean a finite union of triangles (cf. \cite[p.~5]{H}). We
denote by $\pp$ the set of all polygons. A
polygon is {\it simple} if its interior is connected. Two polygons
are {\it nonoverlapping}, if their interiors are disjoint. Every polygon is the
union of finitely many nonoverlapping simple polygons. 

By the {\it direction} of a nonzero vector $v\in \sik$, we mean the unit vector
$v/|v|$. Let ${\cal D}$ be a set of directions in the plane. We denote by
$\pp _{\cal D}$ the family of all polygons $P$ such that the direction of each
side of $P$ belongs to ${\cal D}$. Note that $\pp _{\cal D}$ is a ring in the
sense that if $A,B\in \pp _{\cal D}$, then $A\cup B\in \pp _{\cal D}$, and the
closure of the interior of both $A\cap B$ and of $A\se B$ also belong to
$\pp _{\cal D}$.

We say that the polygons $A,B\in \pp _{\cal D}$ are
{\it ${\cal D}$-equidecomposable} if $A$ can be decomposed into nonoverlapping
polygons $A_1 \stb A_n \in \pp _{\cal D}$ such that suitable translated copies of
$A_1 \stb A_n \in \pp _{\cal D}$ form a decomposition of $B$ into nonoverlapping
polygons. We denote this fact by $A\sim _{\cal D} B$. It is easy to check that
the relation $\sim _{\cal D}$ is an equivalence relation. Our aim is to find
conditions implying $A\sim _{\cal D} B$.

Note that Dehn's theorem is the special case, when ${\cal D}$ only consists of
the direction of the $x$-axis and the $y$-axis. In this case $\pp _{\cal D}$
equals the set of polygons only having sides parallel to the axes. We denote
this set of polygons by $\hh$. Clearly, $A\in \hh$ if and only if $A$ is the
union of finitely many nonoverlapping rectangles with sides parallel to
the axes.

The other extremal case is when ${\cal D}$ is the set of all directions, when
the theorem of Hadwiger and Glur gives the necessary and sufficient condition:
the polygons $A$ and $B$ of equal area are equidecomposable using translations
if and only if $\nu _u (A)=\nu _u (B)$ for every unit vector $u$ (\cite{HG},
\cite[p.~78]{B}). Here $\nu _u$ is an invariant to be defined shortly.

We say that the map $\mu \colon \pp _{\cal D} \to \rr$ is {\it additive} if,
whenever $A\in \pp _{\cal D}$ is decomposed into nonoverlapping polygons
$A_1 \stb A_n \in \pp _{\cal D}$, then $\mu (A)=\sumin \mu (A_i )$. By an
{\it invariant} on $\pp _{\cal D}$ we mean a translation invariant additive
function defined on $\pp _{\cal D}$. (Our notion of invariant coincides with the
notion of translation invariant simple valuation. See \cite{A}, \cite{KP},
\cite{M}, \cite{S}.) Invariants with values in Abelian groups other than $\rr$
can be useful in some contexts (see, e.g. \cite{L}). In these notes, however,
we only consider real valued invariants.

It is clear that if $A,B\in \pp _{\cal D}$ and $A\sim _{\cal D} B$, then
$\mu (A)=\mu (B)$ for every invariant $\mu$. We shall prove the converse:
\begin{theorem}\label{t1}
For every set of directions ${\cal D}$, and for every $A,B\in \pp _{\cal D}$ we
have $A\sim _{\cal D} B$ if and only if $\mu (A)=\mu (B)$ whenever $\mu$ is a
invariant on $\pp _{\cal D}$.
\end{theorem}
This gives the following corollary:
\begin{cor}\label{c1}
\begin{enumerate}[{\rm (i)}]
\item The relation $\sim _{\cal D}$ satisfies the cancellation law. That is, if
$A$ is dissected into the polygons $A_1 \stb A_n \in \pp _{\cal D}$ and $B$ is
dissected into the polygons $B_1 \stb B_n \in \pp _{\cal D}$ such that
$A_1 \sim _{\cal D} A_2 \sim _{\cal D} \ldots \sim _{\cal D} A_n$, $B_1 \sim _{\cal D}
B_2 \sim _{\cal D} \ldots \sim _{\cal D} B_n$ and $A\sim _{\cal D} B$, then
$A_1 \sim _{\cal D} B_1$.
\item  The relation $\sim _{\cal D}$ satisfies the subtraction law. That is, if
$A ,B_1 \in \pp _{\cal D}$ are nonoverlapping, $C ,B_2 \in \pp _{\cal D}$ are
nonoverlapping, and $A\cup B_1 \sim _{\cal D} C\cup B_2$ and $B_1 \sim _{\cal D}
B_2$, then $A\sim _{\cal D} C$. 
\end{enumerate}
\end{cor}
The analogous statements concerning all polytops and an arbitrary group of
isometries containing translations are well-known; see
\cite[Satz VIII, p.~58]{H}, \cite[Theorem 5.1]{KP}. (Note, however, that in
the ring $\hh$ the statement of Theorem \ref{t1} and (i) of Corollary \ref{c1}
may fail for a suitable subgroup of the group of translations. See Example 1
in \cite{L}.) Hadwiger's general theorem is proved by establishing first the
analogue of Corollary \ref{c1}, then by embedding the type semigroup into a
linear space over the rationals, and then using the linear maps of this linear
space to construct the invariants. We follow a partially different route,
mainly because we want to prove a stronger statement involving only some special
invariants. We also want to describe all invariants of
$\pp_{\cal D}$.

Let $L\colon \sik \to \sik$ be a nonsingular linear transformation. If $d$ is a
direction and $v$ is vector of direction $d$, then we denote by $\tilde L (d)$
the direction of the vector $L(v)$; that is, $L(v)/|L(v)|$. It is clear that
$\tilde L (d)$ is well-defined; that is, independent of the choice of $v$. If
${\cal D}$ is a set of directions, then we put
$\tilde L ({\cal D})=\{ \tilde L (d)\colon d\in {\cal D}\}$.

Obviously, if $A\in \pp _{\cal D}$, then $L(A)\in \pp _{\tilde L ({\cal D})}$. It is
also clear that if the polygons $A,B$ are nonoverlapping, then so are $L(A)$
and $L(B)$. Putting these facts together, we obtain the following.
\begin{proposition}\label{p1}
For every set of directions ${\cal D}$ and for every nonsingular linear
transformation $L$ we have $A\sim_{\cal D} B \iff L(A)\sim_{\tilde L ({\cal D})} L(B)$
for every $A,B\in \pp _{\cal D}$.  
\end{proposition}
Dealing with ${\cal D}$-equidecomposability of polygons, the set ${\cal D}$
must contain at least two directions, since every polygon has at least two
sides of different directions.

Suppose ${\cal D}$ only contains two directions. For a suitable nonsingular
linear transformation $L$, $\tilde L ({\cal D})$ consists of the horizontal and
the vertical directions. In this case we have $\pp _{\tilde L ({\cal D})} =\hh$,
the context of Dehn's theorem.

If ${\cal D}$ contains at least three directions, then there is a linear
transformation $L$ such that $\tilde L ({\cal D})$ contains the horizontal and
vertical directions, and the direction of the diagonal $\{ (x,x)\colon
x\in \rr \}$.

Let $S$ denote the set of slopes of the nonvertical directions belonging to
${\cal D}$. Using Proposition \ref{p1}, we can see that in order to prove
Theorem \ref{t1} we may assume that either $S=\{ 0\}$, or $0,1\in S$.

Now let a set $S\su \rr$ be given such that either $S=\{ 0\}$, or $0,1\in S$.
Let ${\cal D}_S$ denote the set of directions containing the vertical direction
and the directions having slopes belonging to $S$. For the sake of brevity, we
write $\pp _S$ and $\sim _S$ instead of $\pp _{{\cal D}_S}$ and $\sim _{{\cal D}_S}$.
{\it Note that} (i) {\it ${\cal D}_S$ always contains the vertical direction by
assumption, and} (ii) {\it $\hh \su \pp _S$ for every $S$.}

Now we construct the invariants we will need in the sequel.
Let $u\in \sik$ be a unit vector. If $[a,b]$ is an oriented segment, then we put
$$\nu _u ([a,b])=
\begin{cases}
0 &\text{if $(b-a)/|b-a| \ne \pm u$},\\
|b-a| &\text{if $(b-a)/|b-a| = u$},\\
-|b-a| &\text{if $(b-a)/|b-a| = -u$}.\\
\end{cases}$$
Let $A$ be a simple polygon, and let the vertices of $A$ be $v_i$ $(i=1\stb k)$
listed counterclockwise. Then we define
\begin{equation}\label{ev}
\nu _u (A)= \sumik \nu _u ([v_{i-1},v_i ]),
\end{equation}
where $v_0 =v_k$.
(The function $\nu _u$ is usually defined using the outer normal $u^\perp$
of $A$ instead of $u$; see, e.g. \cite[Section 2]{KP} or \cite[Section 6]{S}.)

If $A$ is an arbitrary polygon and $A$ is the union of the
nonoverlapping simple polygons $A_1 \stb A_n$, then we define $\nu _u (A)=
\sumin \nu _u (A_i )$. It is clear that $\nu _u$ is a continuous $1$-homogeneous
translation invariant valuation for every unit vector $u$. (See \cite[pp. 79-80]{B}, \cite{KP}, \cite{M} and \cite{S}.) We call the functions $\nu _u$
{\it invariants of the first kind}. 

We also consider another set of invariants. First suppose $S=\{ 0\}$; that is,
$\pp _S =\hh$. We say that $\mu$ is an {\it invariant of the second kind
on $\hh$}, if there are additive functions\footnote{A function $f\colon \rr \to \rr$ is {\it additive} if $f(x+y)=f(x)+f(y)$ for every $x,y\in \rr$.} $f,g\colon \rr \to \rr$ such that
$\mu ([a,b]\times [c,d])=f(b-a)\cd g(d-c)$ for every $a<b$ and $c<d$. (See \cite{E}.)

Next suppose $0,1\in S$. Let $f$ be an additive function $f\colon \rr \to \rr$
such that $f(sx)=s\cd f(x)$ for every $x\in \rr$ and $s\in S$. If $A$ is a
simple polygon with vertices $(x_1 ,y_1 ) \stb (x_n ,y_n )=(x_0 ,y_0 )$ listed
counterclockwise, then we put
\begin{equation}\label{ev1}
\mu _f (A)=\sumin f(x_{i-1}+x_i )\cd f(y_i -y_{i-1}).
\end{equation}
If $A$ is an arbitrary polygon and $A$ is the union of the
nonoverlapping simple polygons $A_1 \stb A_n$, then we define $\mu _f (A)=
\sumin \mu _f (A_i )$. In the next section we prove that $\mu _f$, restricted
to $\pp _S$, is an invariant on $\pp _S$ (see Corollary \ref{c2}).
We call the functions $\mu _f$ {\it invariants of the second kind}. 

\begin{remarks}\label{rs1}
{\rm
1. It is clear that if $\mu _f$ is an invariant of the second kind, then
$\mu_f (tA)=t^2 \cd \mu _f (A)$ for every polygon $A$ and for every positive
rational $t$, where $tA=\{ tx\colon x\in A\}$. That is, invariants of the
second kind are rationally $2$-homogeneous. However, as we will see in the
next section, they only form a proper subset of the rationally $2$-homogeneous
invariants of $\pp _S$. 

2. We also note that if the additive function $f$ is not continuous, then
$\mu _f$ is not an invariant on the set $\pp$ of all polygons, because it is not
additive on $\pp$ (see Proposition \ref{p3}). If, on the other hand, $f$ is
continuous, then it is linear; that is, of the form $c\cd x$ with a constant
$c$ (see \cite[Theorem 2, p. 277]{K}). We show that if $f$ is linear, then
$\mu _f$ is a constant multiple of the area.
In fact, if $f(x)= x/\sqrt 2$ $(x\in \rr )$, then the sum in \eqref{ev1} gives
the area of the simple polygon $A$. Indeed, assuming $x_i >0$ $(i=1\stb n)$,
$(x_{i-1}+x_i )\cd |y_i -y_{i-1}|/2$ equals the area of the trapezoid of vertices
$(0,y_{i-1})$, $(x_{i-1},y_{i-1})$, $(x_i ,y_i ), (0,y_i )$, and it is easy to
check that the sum of these areas with suitable signs gives the area of $A$.

Therefore, an invariant of the second kind on $\pp _S$ is an invariant on $\pp$
if and only if it is a constant multiple of the area.
}
\end{remarks}
Now we can state our main result.
\begin{theorem}\label{t4}
\begin{enumerate}[{\rm (i)}]
\item If $S=\{ 0\}$, then for every $A,B\in \hh$ we have $A\sim _S B$ if and
only if $\mu (A)=\mu (B)$ whenever $\mu$ is an invariant of 
the second kind on $\hh$.
\item Suppose that $0,1\in S$, and let $K$ denote the field generated by
$S$. Then for every $A,B\in \pp _S$ we have $A\sim _S B$ if and only if
$A\sim _K B$ if and only if $\mu (A)=\mu (B)$ whenever $\mu$ is an invariant
of the first or of the second kind on $\pp _K$.
\end{enumerate}
\end{theorem}

\begin{remarks}\label{rs2}
{\rm 1.
If $S_1 ,S_2$ are different subsets of $\rr$, then the equivalence
relations $\sim _{S_1} , \sim_{S_2}$ are different, as their domains, $\pp _{S_1} ,
\pp_{S_2}$ are different. However, it may happen that $S_1 \subsetneq S_2$, but
$A\sim _{S_1}B\iff A \sim_{S_2}B$ for every $A,B\in \pp _{S_1}$. This means that
whenever $A,B\in \pp _{S_1}$ are equidecomposable using directions from
${\cal D}_{S_2}$ then they are also equidecomposable using directions from the
smaller set ${\cal D}_{S_1}$. According to the
theorem above, this is always the case if $0,1\in S_1$ and $S_2$ is the field
generated by $S_1$. 

2. Theorem \ref{t4} is sharp. If $S\su S'$ and $A\sim _{S} B\iff A \sim_{S'} B$
for every $A,B\in \pp _S$, then necessarily $S'$ is contained in the field $K$
generated by $S$. See Corollary \ref{c3}.

3. Theorem \ref{t4} is a generalization of the Hadwiger-Glur theorem \cite{HG}.
Indeed, if $S=\rr$ then, as we noted in Remark \ref{rs1}.2, the only invariants
of the second kind are the constant multiples of the area. Therefore,
$A\sim _{\rr} B$ if and only if $A, B$ are of the same
area, and  $\nu _u (A) =\nu _u (B)$ for every unit vector $u$, which is the
theorem of Hadwiger and Glur.

4. In (ii) of Theorem \ref{t4}, the condition $0,1\in S$ cannot be omitted;
that is, statement (ii) of the theorem is not true if $S=\{ 0\}$. Let $s$ be an
irrational number, and put 
\begin{equation*}
R_1 =[0,s ]\times \nl , \ R_2 =\nl \times [0,s ] .
\end{equation*}
If $S=\{ 0\}$, then, by Dehn's theorem, $R_1 \sim_{S} R_2$ is not true, since
$s$ is irrational. However, we have $R_1 \sim_{\qq} R_2$ (see Theorem \ref{t5}
below).

}
\end{remarks}

If $A,B$ are rectangles, then the condition formulated in Theorem \ref{t4}
can be made more explicit. The following theorem is a generalization of
\cite[Satz XIII, p.~76]{H} in dimension two.
\begin{theorem}\label{t5}
Suppose $0,1\in S$. Then the rectangles $R_i =[0,a_i ]\times [0,b_i ]$
$(i=1,2)$ are $S$-equidecom\-posable if and only if $a_1 b_1 =a_2 b_2$,
and at least one of the numbers $a_2 /a_1 \ ( =b_1 /b_2 )$ and 
$b_2 /a_1 \ (= a_2 /b_1 )$ belongs to the subfield of $\rr$ generated by
$S$. 
\end{theorem}
\begin{remark}\label{r2}
{\rm Note that the condition $0,1\in S$ is essential. If $S=\{ 0\}$
then, by Dehn's theorem, $R_i =[0,a_i ]\times [0,b_i ]$ $(i=1,2)$ are
$S$-equidecom\-posable if and only if $a_1 b_1 =a_2 b_2$ and
$a_2 /a_1 \in \qq$. The condition formulated in Theorem \ref{t5} is weaker,
and is not sufficient if $S=\{ 0\}$.
}
\end{remark}

As an application of Theorem \ref{t5} we also prove the following.
\begin{theorem}\label{t6}
Suppose $0,1\in S$, and let $K$ denote the field generated by $S$. Then
\begin{enumerate}[{\rm (i)}]
\item the rectangle $[0,a]\times [0,b]$ is $S$-equidecomposable to a square
with sides parallel to the axes if and only if $a/b=\de ^2$, where $\de \in K$;
\item the rectangle $[0,a]\times [0,b]$ is $S$-equidecomposable to a square
(of arbitrary position) if and only if $a/b=\ga ^2 +\de ^2$, where $\ga ,\de
\in K$.
\end{enumerate}
\end{theorem}
For example, if $S=\qq$, then $\nl \times [0,2]$ is $S$-equidecomposable to a
square, but is not $S$-equidecomposable to a square with sides parallel to the
axes. On the other hand, $\nl \times [0,3]$ is not $S$-equidecomposable to any
square. Furthermore, if $S=\qq (\sqrt 2 )$, then $\nl \times [0,3]$ is $S$-equidecomposable to a
square, but is not $S$-equidecomposable to a square with sides parallel to the
axes. If $S=\qq (\sqrt 3 )$, then $\nl \times [0,3]$ is $S$-equidecomposable
to a square with sides parallel to the axes. 

The structure of the paper is the following. In the next section we construct
rationally $2$-homogeneous invariants (including the invariants of the second
kind) for every $\pp _S$. In Section \ref{s2}
we prove Theorem \ref{t5}. We prove Theorem \ref{t4} in Section \ref{s4}, using
two lemmas presented in Section \ref{s3}. Theorem \ref{t6} will be proved in
Section \ref{s5}. The description of all invariants on $\pp_S$, through an
analogue of the Hadwiger-McMullen decomposition, is given in
Section \ref{s6}. In the last section we make some comments on the duality
between the rings $\pp _S$ (that is, subsets of $\rr$) and the sets of
invariants (that is, sets of symmetric biadditive functions).

\section{Construction of the rationally $2$-homogeneous invariants} \label{s1}
In this section we construct rationally $2$-homogeneous invariants on
$\pp _S$ in the case when $0,1\in S$. As we will prove in Corollary \ref{c4},
every rationally $2$-homogeneous invariant on $\pp _S$ can be
obtained this way.

Let $F\colon \sik \to \rr$ be a biadditive function. If $[a,b]$ is an
oriented segment, where $a=(a_1 ,a_2 )$ and $b=(b_1 ,b_2 )$, then we define
\begin{equation*}
\mu _F ([a,b])= F(a_1 +b_1 ,b_2 -a_2 ). 
\end{equation*}
Let $A$ be a simple polygon. Then we define
\begin{equation*}
\mu _F (A)= \sumik \mu _F ([v_{i-1},v_i ]),
\end{equation*}
where $v_1 \stb v_k =v_0$ are the vertices of $A$ listed counterclockwise.
If $A$ is an arbitrary polygon and $A$ is the union of the nonoverlapping
simple polygons $A_1 \stb A_n$, then we put $\mu _F (A)=\sumin \mu _F (A_i )$.
In this way we have defined a map $\mu _F \colon \pp \to \rr$ from the
family $\pp$ of all polygons into $\rr$.
\begin{proposition}\label{p2}
The function $\mu _F$ is invariant under translations.
\end{proposition}
\proof It is clear that
\begin{align*}
  \mu _F ([a,b] +c)= & F(a_1 +b_1 +2c_1 ,b_2 -a_2 )=\\
  &F(a_1 +b_1 ,b_2 -a_2 )+ F(2c_1 ,b_2 -a_2 )=\\
  &\mu _F ([a,b]) + F( 2c_1 ,b_2 -a_2 )
\end{align*}
for every $a=(a_1 ,a_2 ) ,b=(b_1 ,b_2 ),c=(c_1 ,c_2 )\in \sik$. If $A$ is a simple polygon with vertices
$v_i =(x_i ,y_i )$ $(i=1\stb k)$ listed counterclockwise, then we obtain
\begin{align*}
  \mu _F (A+c)=&\sumik \mu _F ([v_{i-1},v_i ] +c)=\\
  & \sumik \mu _F ([v_{i-1},v_i ]) +\sumik F(2c_1 ,y_i -y_{i-1}) =\\
  & \mu _F (A)+F\left( 2c_1 ,\sumik (y_i -y_{i-1}) \right)=\\
  & \mu _F (A)+F( 2c_1 ,0)=\mu _F (A) .
\end{align*}
Therefore, we have $\mu _F (A+c)=\mu _F (A)$ for every
$A\in \pp$. \hfill $\square$

\begin{theorem}\label{t7}
Suppose that $0,1\in S$. Let $F\colon \sik \to \rr$ be a symmetric,
biadditive function such that
\begin{equation}\label{e1}
F(x,sy)=F(sx,y) \qquad (s\in S, \ x,y\in \rr ).
\end{equation}
Then $\mu _F$ is an invariant on $\pp _S$. 
\end{theorem}
\proof We only have to show that $\mu _F$ is additive on $\pp _S$. First we
prove that if the segment $[a,b]$ is vertical or its slope belongs to $S$, then 
\begin{equation}\label{e3}
\mu _F ([a,b])= \mu _F ([a,c])+\mu _F ([c,b]) 
\end{equation}
for every $c\in [a,b]$. Let $a=(a_1 ,a_2 )$, $b=(b_1 ,b_2 )$ and $c=(c_1 ,c_2 )$.
If the segment $[a,b]$ is vertical, then $a_1 =c_1 =b_1$, and \eqref{e3}
follows from $F(2a_1 ,b_2 -a_2 )= F(2a_1 ,c_2 -a_2 )+ F(2a_1 ,b_2 -c_2 )$, which
is clear from the biadditivity of $F$. Suppose $[a,b]$ is not vertical,
and let its slope be $s\in S$. Then $c_2 -a_2 =s(c_1 -a_1 )$, and thus
\begin{align*}
  \mu _F ([a,c])=& F(a_1 +c_1 ,s(c_1 -a_1 ))=\\
  &F(a_1 , sc_1 )-
  F(a_1 , sa_1 )+F(c_1 , sc_1 )-F(c_1 , sa_1 )=\\
&  F(c_1 , sc_1 )-F(a_1 , sa_1 ),
\end{align*}
since $F(a_1 , sc_1 )=F(sc_1 ,a_1 )=F(c_1 , sa_1 )$. We obtain
$\mu _F ([c,b])=F(b_1 , sb_1 )-F(c_1 , sc_1 )$ and $\mu _F ([a,b])=F(b_1 , sb_1 )-F(a_1 , sa_1 )$ the same way. Therefore,
\begin{align*}
  \mu _F ([a,c])&+ \mu _F ([c,b]) =\\
  & (F(c_1 ,sc_1 ) -F(a_1 ,sa_1 )) +(F(b_1 ,sb_1 ) -F(c_1 ,sc_1 ))=\\
  & F(b_1 ,sb_1 ) -F(a_1 ,sa_1 )=\mu _F ([a,b]).
\end{align*}
This proves \eqref{e3}. Now let $A=\sumin A_i$ be a decomposition of
$A\in \pp _S$ into the nonoverlapping polygons $A_i$ belonging to $P_S$.
We have to prove
\begin{equation}\label{e4}
\mu _F (A)= \sumin \mu _F (A_i ).
\end{equation}
By the definition of $\mu _F$, $\sumin \mu _F (A_i )$ is a sum of the
form $\sum \mu _F ([a,b])$, where $[a,b]$ runs through the sides of the sets
$A_i$. The orientation of the segments $[a,b]$ is obtained by listing the
vertices of the corresponding $A_i$ counterclockwise.

Let $V$ denote the set of vertices of the polygons $A_i$ $(i=1\stb n)$. 
It follows from \eqref{e3} that placing extra points on the sides of the
polygons $A_i$ and treating them as vertices, the value of the sum
$\sum \mu _F ([a,b])$ does not change. Therefore, placing the points of $V$ on
every side of $A_1 \stb A_n$ that contains them, the value of the sum
$\sum \mu _F ([a,b])$ does not change. In the new sum $[a,b]$ runs
through all segments lying on the union of the boundaries of $A_1 \stb A_n$
and such that $a,b\in V$. If such a segment $[a,b]$ is on the boundary of
one of the sets $A_i$ and lies in the interior of $A$ except perhaps its
endpoints, then the segment $[b,a]$ also appears in the sum, as part of the
boundary of another polygon $A_j$. Since $\mu _F ([a,b])=-\mu _F ([b,a])$, 
the sum of these terms is zero. Thus $\sum \mu _F ([a,b])$ equals the sum
of those terms $\mu _F ([a,b])$ for which $[a,b]$ lies on the boundary of $A$.
The value of this sum is $\mu _F (A)$, proving \eqref{e4}.  \hfill $\square$

\begin{cor}\label{c2}
Let $S\su \rr$ be given with $0,1\in S$, and let $f\colon \rr \to \rr$ be an
additive function such that $f(sx)=s\cd f(x)$ for every $x\in \rr$ and $s\in S$.
Then $\mu _f$ is an invariant on $\pp _S$.
\end{cor}
\proof
Put $F(x,y)= f(x )\cd f(y)$ for every $x,y\in \rr$. It is clear that $F$ is a
symmetric biadditive function satisfying \eqref{e1}, and then we can apply
Theorem \ref{t7}. \hfill $\square$

\begin{remark}\label{r4}
{\rm What we call invariants of the second kind only constitutes a proper subset
of all invariants constructed in Theorem \ref{t7}. Indeed, symmetric
biadditive functions $F(x,y)$ are not necessarily of the form $f(x)\cd f(y)$.
}
\end{remark}

\begin{proposition}\label{p3}
Let $f\colon \rr \to \rr$ be additive. Then $\mu _f$ is additive on the set
$\pp$ of all polygons if and only if $f$ is linear. 
\end{proposition}

\proof We saw already in the introduction that if $f$ is linear, then $\mu _f$
is a constant multiple of the area (see Remark \ref{rs1}.2) and, consequently,
$\mu _f$ is an invariant on $\pp$.

Therefore, it is enough to show the converse. It is easy to check
that if $A=[0,x]\times [0,x]$ and $B=\nl \times [0,x^2 ]$, where $x>0$, then
$\mu _f (A)=2f(x)^2$ and $\mu _f (B)=2f(1)f(x^2 )$. If $\mu _f$ is additive on
$\pp$ then, as $A$ and $B$ are equidecomposable in $\pp$, we have
$f(x)^2 =f(1)f(x^2 )$ for every $x>0$.
Since $f$ is odd, this is true for every $x\in \rr$.

If $f(1)=0$, then we obtain $f\equiv 0$, and so $f$ is linear in this case.
If $f(1)\ne 0$, then $f(x^2 )=f(x)^2 /f(1)$ for every $x$, and thus either
$f\ge 0$ or $f\le 0$ on $[0,\infty )$. Since $f$ is additive, this implies
that $f$ is linear (see \cite[Theorem 2, p. 277]{K}).  \hfill $\square$

\section{Rectangles}\label{s2}
In this section we prove Theorem \ref{t5}.

{\bf I.} First we prove the `only if' statement of the theorem.

Let the rectangles $R_i =[0,a_i ]\times [0,b_i ]$ $(i=1,2)$ be
$S$-equidecom\-posable. Then the area of $R_1$ equals that of $R_2$; that is,
$a_1 b_1 =a_2 b_2$. We have to prove that at least one of the
numbers $a_2 /a_1 \ (=b_1 /b_2 )$ and $b_2 /a_1  \ (=b_1 /a_2 )$ belongs to the
field generated by $S$.

Let $f$ be an additive function satisfying $f (sx)=s\cd f (x)$ for every
$s\in S$ and $x\in \rr$. By Corollary \ref{c2}, $\mu _f$ is an invariant
on $\pp _S$.

It is easy to check that $\mu _f (R_i )=2f(a_i )\cd f(b_i )$ $(i=1,2)$. If
$R_1 \sim _S R_2$, then the invariance of $\mu _f$ implies
\begin{equation}\label{e5}
f (a_1 )\cd f (b_1 )=f (a_2 )\cd f (b_2 ).
\end{equation}
Let $K$ denote the field generated by $S$, and suppose that $a_2 /a_1 \notin K$
and $b_2 /a_1 \notin K$. 
First we assume that $a_1 ,a_2$ and $b_1$ are linearly independent over $K$.
Then there is a basis $B$ of $\rr$ as a linear space over the field $K$ such
that $a_1 ,a_2 ,b_1 \in B$. Then every real number $x$ has a unique
representation $\sum_{b\in B} \al _b (x)\cd b$, where $\al _b (x)\in K$ for every
$b\in B$, and $\al _b (x)=0$ for all but a finite number of $b\in B$. We define
$f (x)=\al _{a_1} (x)+\al _{b_1} (x)$ for every $x\in \rr$. It is easy to see
that $f $ is additive, and $f (sx)=s\cd f(x)$ for every $s\in K$ and
$x\in \rr$. Thus \eqref{e5} must hold. However, we have $f(a_1 )=1$,
$f(b_1 )=1$ and $f(a_2 )=0$, which is a contradiction.

Next suppose that $a_1 ,a_2$ and $b_1$ are linearly dependent over $K$.
Since $a_2 /a_1 \notin K$ by assumption, this implies $b_1 =\la _1 a_1 + \la _2
a_2$, where $\la _1 ,\la _2 \in K$. The assumption $b_2 /a_1 =a_2 /b_1 \notin K$
implies $\la _1 \ne 0$. There is a basis $B$ of $\rr$ as a linear space over
the field $K$ such that $a_1 ,a_2  \in B$. Then every real number $x$ has a
unique representation $\sum_{b\in B} \al _b (x)\cd b$ as above. We define
$f (x)=\al _{a_1} (x)$ for every $x\in \rr$. Then $f$ is additive and
$f (sx)=s\cd f(x)$ for every $s\in K$ and $x\in \rr$. Then \eqref{e5}
must hold. However, we have $f(a_1 )=1$, $f(b_1 )=\la _1 \ne 0$ and
$f(a_2 )=0$, which is a contradiction. This proves the `only if' part
of the theorem.

{\bf II.} In order to prove the `if' statement of the theorem we present some lemmas on the type
semigroup. (Type semigroups appeared in Tarski's work \cite{T}; see also
\cite[p.~168]{TW}.)

If $A\in \pp _S$, then we denote by $[A]_S$ the set
$\{ B\in \pp _S \colon B\sim _S A\}$. If $A,B\in \pp _S$, then we define
$[A]_S +[B]_S =[A\cup B']_S$, where $B' \sim _S B$ and $A$ and $B'$ are
nonoverlapping. It is easy to check that $[A]_S +[B]_S $ does not depend on the
choice of $B'$. Then $T_S =\{ [A]_S \colon A\in \pp _S\}$ equipped with this
addition is a commutative semigroup. In the sequel, if $S$ is clear from the
context, then we will write $[A]$ instead of $[A]_S$.
\begin{lemma}\label{l1}
Let $S$ be arbitrary.
\begin{enumerate}[{\rm (i)}]
\item If $H\in \hh$, and $\phi$ is a homothetic transformation of ratio $n$,
where $n$ is a positive integer, then $[\phi (H)]=n^2 \cd [H]$.
\item If $A,B\in \hh$ and there is a positive integer $n$ such that
$n\cd [A]=n\cd [B]$, then $A\sim _S B$.
\end{enumerate}
\end{lemma}
\proof
(i) Let $\phi$ be a homothetic transformation of ratio $n$. If $R$ is
a rectangle, then $\phi (R)$ can be decomposed into $n^2$ congruent rectangles,
each of which is a translated copy of $R$.

Let $H=\bigcup_{i=1}^k A_i$, where $A_1 \stb A_k$ are nonoverlapping
rectangles with sides parallel to the axes. Then $\phi (A_i )$ can be
decomposed into $n^2$ translated copies of $A_i$ for every $i$. This implies
that $[\phi (H)]=n^2 \cd [H]$.

(ii) If $n\cd [A]=n\cd [B]$, then $n^2 \cd [A]=n^2 \cd [B]$.
By (i) this implies $[\phi (A)]=[\phi (B)]$; that is,
$\phi (A)\sim _S \phi (B)$. By Proposition \ref{p1}, we obtain
$A\sim _S B$, as the homothetic transformation $\phi$ does not change the set $S$ of directions. \hfill $\square$

If $\mathfrak G$ is a commutative semigroup and $x,y\in \mathfrak G$, then we write $x\le y$
if there is a $z\in \mathfrak G$ such that $x+z=y$.
\begin{lemma}\label{l2}
Let $\mathfrak G$ be a commutative semigroup, and let $a,b,c\in \mathfrak G$. Suppose that  
$a+b=c+b$ and $b\le na$, $b\le nc$ for a positive integer $n$. Then we have
$ka=kc$, where $k=n^2 +n$.
\end{lemma}
\proof
First note that if $x,y,z\in \mathfrak G$ and $x+z=y+z$, then $ix+z=iy+z$ for every
positive integer $i$. Indeed, this is true for $i=1$, and if it is true for
$i$, then
\begin{align*}
  (i+1)x+z =&x+(ix +z)=x+(iy+z)=iy+(x+z)=\\
  &iy+(y+z)=(i+1)y+z.
\end{align*}
Turning to the proof of the lemma, we have $na=a_1 +b$ and $nc=c_1 +b$ for
some $a_1 ,c_1 \in \mathfrak G$. From $a+b=c+b$ we obtain $na+nb=nc+nb$ and $a_1 +(n+1)b=
c_1 +(n+1)b$. Applying the observation above we get $(n+1)a_1 +(n+1)b=
(n+1)c_1 +(n+1)b$, and thus
\begin{align*}
  (n^2 +n)a=&(n+1)\cd na =(n+1)(a_1 +b)=\\
  &(n+1)(c_1 +b)=\\
&(n+1)\cd nc= (n^2 +n)c. \hskip 5cm \square
\end{align*}
\begin{lemma}\label{l3}
Let $S$ be arbitrary, and let $H_1 , H_2 , C_1 , C_2$ be nonoverlapping
polygons such that $H_1 ,H_2 \in \hh$, $C_1 , C_2 \in \pp _S$, $C_1 \sim _S
C_2$, and $H_1 \cup C_1 \sim _S H_2 \cup C_2$. Then $H_1 \sim _S H_2$.
\end{lemma}
\proof
In the formalism of the type semigroup $T_S$ we have $[C_1 ]=[C_2 ]$
and $[H_1 ]+[C_1 ]=[H_2 ]+ [C_1]$. There are homothetic transformations 
$\phi _1 ,\phi _2$ of ratio $n$ with a suitable positive integer $n$ such that
$C_1 \su \phi_1 (H_1)$ and $C_1 \su \phi_2 (H_2 )$.
This implies $[C_1 ]\le n^2 [H_1 ]$ and $[C_1 ]\le n^2 [H_2 ]$. By Lemma \ref{l2}, we have
$k\cd [H_1 ]=k\cd [H_2 ]$ with a suitable positive integer $k$. Then Lemma
\ref{l1} gives $H_1 \sim _S H_2$. \hfill $\square$

{\bf III.} Now we present two lemmas on the equidecomposability of rectangles.
Recall that under our assumptions we have either $S=\{ 0\}$ or $0,1\in S$.
\begin{lemma}\label{l4}
\begin{enumerate}[{\rm (i)}]
\item For every $S$, we have $([0,a]\times [0,b])\sim _S ([0,ra]\times [0,b/r])$
for every $a,b>0$ and $r\in \qq$, $r>0$.
\item If $0,1\in S$, then
$$([0,a]\times [0,b])\sim _S ([0,b/t]\times [0,t a])
  \sim _S ([0,tb]\times [0,a/t]) $$ 
for every $a,b,t>0$, where $t=|s|$ for some $s\in S$.
\item If $0,1\in S$, then $([0,a]\times [0,b])\sim _S ([0,b]\times [0,a])$
for every $a,b>0$.
\end{enumerate}
\end{lemma}
\proof
(i): If $r=p/q$, where $p,q$ are positive integers, then both $[0,a]\times
[0,b]$ and $[0,ra]\times [0,b/r]$ can be dissected into $p\cd q$ translated
copy of $[0,a/q]\times [0,b/p]$. This implies that $[0,a]\times [0,b]$ and
$[0,ra]\times [0,b/r]$ are equidecomposable (even in $\hh$).

(ii) and (iii): Let $a,b>0$ and $s\in S$, $s\ne 0$ be given. Put $t=|s|$. By
(i), we may replace $[0,a]\times [0,b]$ by $[0,ra]\times [0,b/r]$ for every
positive rational $r$. Choosing an appropriate $r$, we may assume that $b/a<t<
b/(a-b)$. This condition guarantees that Figure \ref{fig1} applies, and shows
a correct proof of
\begin{equation}\label{e5a}
([0,a]\times [0,b]) \sim_S ([0,b/t]\times [0,ta]).
\end{equation}
Since $1\in S$ by assumption, \eqref{e5a} gives (iii). 
Then, applying \eqref{e5a} again we obtain
$$([0,tb]\times [0,a/t]) \sim _S ([0,a/t] \times[0,t b])
\sim _S ([0,b]\times [0,a])\sim _S ([0,a]\times [0,b]),$$
which completes the proof of (ii). \hfill $\square$

\begin{figure}[h!]
\centering
\begin{tikzpicture}[scale=0.5]

\coordinate (A) at (0,0);
\coordinate (B) at (10,0);
\coordinate (C) at (10,11);
\coordinate (D) at (0,11);

\coordinate (C_1) at (4.545,11);
\coordinate (C_2) at (10,6);

\coordinate (D_1) at (4.545,0);
\coordinate (D_2) at (0,6);

\coordinate (E) at (4.545,6);

\coordinate (M) at (2,4);

\coordinate (N) at (6,8);


\draw (A) -- (B) -- (C) -- (C_1)--(E) -- (D_2) -- cycle;
\draw (D_1) -- (E) -- (C_2);
\draw (A)--(C);


\draw[-{Latex[length=3mm]}] (M) -- (N);


\node[left] at (0, 3) {$b$};
\node[below] at (5,0) {$a$};
\node[right] at (10, 5) {$s a$};
\node at (7,11.6){$b/s$}; 
\node at (1.2,13){If $s>0:$}; 


\coordinate (A') at (14,0);
\coordinate (B') at (24,0);
\coordinate (C') at (24,11);
\coordinate (D') at (14,11);

\coordinate (C_1') at (19.454,11);
\coordinate (C_2') at (24,6);

\coordinate (D_1') at (19.454,0);
\coordinate (D_2') at (14,6);

\coordinate (E') at (19.454,6);

\coordinate (M') at (22,4);

\coordinate (N') at (18,8);


\draw (A') -- (B') -- (C_2') -- (E') -- (C_1') -- (D')--cycle;
\draw (D_1') -- (E') -- (D_2');
\draw (B')--(D');



\draw[-{Latex[length=3mm]}] (M') -- (N');


\node[left] at (14, 5) {$|s|a$};
\node[below] at (19,0) {$a$};
\node[right] at (24, 3) {$b$};
\node at (17,11.6){$b/|s|$}; 
\node at (15.5,13){If $s<0:$};
\end{tikzpicture}

\caption{Equidecomposition of  
$[0,a]\times [0,b]$ and $[0,b/|s|]\times
[0,|s| \cd a]$}
    \label{fig1}
\end{figure}

\begin{lemma}\label{l5}
Suppose $0,1\in S$, and let $K$ denote the field generated by $S$. If
$a\in K$, $a>0$, then $([0,a]\times [0,b]) \sim_S (\nl \times [0,ab])$ for every
$b>0$.
\end{lemma}
\proof
Let
$$X=\{ x\in \rr , \ x>0 \colon ([0,x]\times [0,b]) \sim_S ( \nl \times [0,bx])
\ \text{for every} \ b>0 \} .$$
We have to prove that $a\in X$ for every $a\in K$, $a>0$.
It follows from Lemma \ref{l4} that $X$ contains the positive rationals,
and also the elements of the set $S^+ =\{ |s|\colon s\in S, \ s\ne 0\}$.
We prove that $X$ is closed under addition and multiplication.

Let $a,b\in X$. For every $c>1$ we have $([0,a]\times [0,c]) \sim_S (\nl
\times [0,ac])$ and $([0,b]\times [0,c]) \sim_S (\nl \times [0,bc])$. Since
$[0,a+b]\times [0,c]=([0,a]\times [0,c])\cup ([a,a+b] \times [0,c])$ and
$\nl \times [0,(a+b)c] =([0,1]\times [0,ac])\cup (\nl \times [ac,(a+b)c]$,
we get $([0,a+b]\times [0,c])\sim _S (\nl \times [0,(a+b)c])$ and $a+b \in X$.

Now $a\in X$ gives $([0,a]\times [0,c/b]) \sim_S (\nl \times [0,ac/b])$.
Applying a homothetic transformation of ratio $b$ we obtain
$([0,ab]\times [0,c]) \sim_S ([0,b] \times [0,ac])$. By $b\in X$ we get
$([0,b]\times [0,ac]) \sim_S (\nl \times [0,abc])$. Thus $([0,ab]\times [0,c])
\sim_S (\nl \times [0,abc])$, and thus $ab\in X$.

Let $\qq [S]^+$ denote the set of numbers of the form $\sumin r_i \cd t_i$,
where $r_i$ is a positive rational, and $t_i$ is a product of elements of
$S^+$ for every $i$. Clearly, every element of the ring $\qq [S]$ is the
difference of two elements of $\qq [S]^+$. Since $X$ is closed under addition
and multiplication, we have $\qq [S]^+ \su X$. Next we prove that if
$a\in \qq [S]$ and $a>0$, then $a\in X$.

Let $a=a_1 -a_2$, where $a_1 ,a_2 \in \qq [S]^+$. Let $b>0$ be given, and put

\begin{tabular}{ll}
   $A=[0,a]\times [0,b]$,  &$C=\nl \times [0,ab]$ , \\
   $B_1 =[a,a+a_2 ]\times [0,b]$, &$B_2 =\nl \times [ab,ab+ a_2 b]$.
\end{tabular}

Since $a+a_2 =a_1 \in \qq [S]^+ \su X$, we have
$$A\cup B_1 =([0,a_1 ]\times [0,b]) \sim _S ([0,1]\times [0,a_1 b])=C\cup B_2$$
Also, $B_1 \sim _S B_2$ by $a_2 \in X$. Then, by Lemma \ref{l3}, we get
$A\sim _S C$ and $a\in X$.

Now we can complete the proof of the lemma. Let $a\in K=\qq (S)$, $a>0$ be
given, and let $b>0$ be arbitrary. Then $a=a_3 /a_4$, where $a_3 ,a_4$ are
positive elements of $\qq [S]$. Then we have $a_3 ,a_4 \in X$, and thus
\begin{align*}
  & ([0,a_3 ]\times [0,a_4 b]) \sim _S (\nl \times [0,a_3 a_4 b]) ,\\
  & ([0,a_4 ]\times [0,a_3 b]) \sim _S (\nl \times [0,a_3 a_4 b]) ,\\
\end{align*}
hence $([0,a_3 ]\times [0,a_4 b]) \sim _S ([0,a_4 ] \times [0,a_3 b])$.
Applying a homothetic transformation of ratio $1/a_4$ we obtain
$([0,a_3 /a_4 ]\times [0,b]) \sim _S (\nl \times [0,a b])$. Thus $a\in X$,
and the proof is complete. \hfill $\square$

{\bf IV.}
We turn to the proof of the `if' part of Theorem \ref{t5}. Let $S\su \rr$ be
such that $0,1\in S$. Suppose that the rectangles
$R_i =[0,a_i ]\times [0,b_i ]$ $(i=1,2)$ are such that $a_1 b_1 =a_2 b_2$,
and either $a_2 /a_1 \in K$ or $b_2 /a_1 \in K$, where $K$ is the field
generated by $S$. We have to prove that $R_1 \sim _S R_2$.

By (iii) of Lemma \ref{l4} we have
$([0,a_1 ]\times [0,b_1 ]) \sim _S ([0,b_1 ]\times [0,a_1 ])$. Therefore, we may
replace $R_1$ by $[0,b_1 ]\times [0,a_1 ]$ if necessary, and we may assume
that $a_2 /a_1 \in K$. If $L(x,y)=(x/a_2 ,y)$ $(x,y\in \rr )$, then
$L$ is a linear transformation, and $L(R_1 )=[0,a_1 /a_2 ]\times [0,b_1 ]$
and $L(R_2 )=[0,1]\times [0,b_2 ]$. Since $a_1 /a_2 \in K$ and $b_2 =
a_1 b_1 /a_2$, (ii) of Lemma \ref{l4} gives
$$L(R_1 )\sim _S (\nl \times [0,a_1 b_1 /a_2 ])=(\nl \times [0,b_2 ]) =L(R_2 ).$$
Then, by Proposition \ref{p1}, we obtain $R_1 \sim _S R_2$. \hfill $\square$

\section{Further lemmas}\label{s3}
\begin{lemma}\label{l6}
Let $K$ be a subfield of $\rr$, and let $a_1 \stb a_n$ be positive real
numbers. Then there are positive real numbers $b_1 \stb b_t$ that are linearly
independent over $K$ and such that for every $i=1\stb n$, $a_i$ is a linear
combination of $b_1 \stb b_t$ with nonnegative coefficients belonging to $K$.
\end{lemma}
\proof  
If $b_1 \stb b_t \in \rr$, then we denote by $L_+ (b_1 \stb b_t )$ the set of
numbers $\sum_{i=1}^t r_i \cd b_i$, where $r_i \ge 0$ and $r_i \in K$
$(i=1\stb t)$.

We start with three simple observations. First note that multiplying the
numbers $a_1 \stb a_n$ by positive numbers belonging to $K$ affects neither the condition nor the conclusion of the statement of the lemma.

Next note that if $a_1 \stb a_n$ are linearly independent over $K$, then we
can take $t=n$ and $b_i =a_i$ $(i=1\stb n)$.

The third fact is the following: {\it if $c, c_1 \stb c_k$ are positive
real numbers such that $c<\sumik c_i$, then there are
positive rational numbers $s_i$ such that $\sumik s_i =1$ and
$s_i \cd c<c_i$ $(i=1\stb k)$.} Indeed, let $\ep >0$ be such that
$\sumik (c_i -\ep )>c$, and choose positive rational numbers $r_i$
such that $c_i -\ep < r_i \cd c<c_i$ $(i=1\stb k)$. Then $\sumik r_i \cd c>
\sumik (c_i -\ep )>c$, $\sumik r_i >1$, and thus the numbers $s_j =
r_j /\sumik r_i$ $(j=1\stb k)$ will satisfy the requirements.

We turn to the proof of the lemma. We prove by induction on $n$.
The case $n=1$ is obvious.

Let $n>1$, and suppose that the statement is true for $n-1$ positive numbers.
If $a_1 \stb a_n$ are linearly independent over $K$, then we are done. Suppose
this is not the case, and let $\sumin r_i \cd a_i =0$, where $r_1 \stb r_n
\in K$ and not all $r_i$ are zero. Since the numbers $a_i$ are positive, some
of the coefficients $r_1 \stb r_n$ are positive, and some of them are negative.
We may assume that $r_i >0$ if $i=1\stb k$, $r_i <0$ if $i=k+1 \stb m$, and
$r_i =0$ if $m+1 \le i\le n$, where $0<k<m\le n$. Replacing $a_i$ by
$a_i /|r_i |$ for every $i=1\stb m$ we may also assume that $r_i =1$
for every $i=1\stb k$, and $r_i =-1$ for every $i=k+1\stb m$. That is, we have
\begin{equation}\label{e6}
a_1 +\ldots +a_k =a_{k+1} +\ldots + a_{m} ,
\end{equation}
where $0<k<m$. We prove the statement by induction on $m$.
If $m=k+1$ then take the system of numbers $\{ a_i \colon 1\le i\le n, \
i\ne k+1\}$. It contains $n-1$ elements and thus, by the induction hypothesis
on $n$, we obtain positive numbers $b_1 \stb b_t$ such that they are
linearly independent over $K$, and $a_i \in L_+ (b_1 \stb b_t )$ for every
$1\le i\le n, \ i\ne k+1$. By \eqref{e6}, we also have
$a_{k+1} \in L_+ (b_1 \stb b_t )$, and we are done.

Next suppose that $m>k+1$ and that the statement is true when \eqref{e6}
holds with $m-1$ in place of $m$. Since $a_m <a_1 +\ldots + a_k$, there
are positive rational numbers $s_i$ $(i=1\stb k)$ such that
$\sumik s_i=1$ and $s_i \cd a_m < a_i$ $(i=1\stb k)$.
Then we have
\begin{equation}\label{e6a}
\sumik (a_i -s_i \cd a_m )=a_{k+1} +\ldots + a_{m-1} .
\end{equation}
Now take the system of numbers
$$Z= \{ a_i -s_i \cd a_m \colon 1\le i\le k\} \cup \{ a_i \colon k+1 \le i\le n
\} .$$
It contains $n$ elements. By \eqref{e6a}, and by the induction hypothesis on
$m$, we obtain positive numbers $b_1 \stb b_t$ such that they are linearly
independent over $K$, and $Z\su L_+ (b_1 \stb b_t )$.  
Since $a_i =(a_i -s_i \cd a_m )+s_i \cd a_m$ for every $i=1\stb k$, we have
$a_i \in L_+ (b_1 \stb b_t )$ for every $i=1\stb n$. \hfill $\square$

Recall that $\hh$ denotes the set of polygons only having sides parallel
to the axes. The invariants $\nu _u$ were defined in \eqref{ev}.

\begin{lemma}\label{l7a}
Suppose $0,1\in S$. For every $A\in \pp _S$ there are nonoverlapping
polygons $H, T_1 \stb T_k \in \pp _S$ such that $A\sim _S H\cup T_1 \cup
\ldots \cup T_k$, $H\in \hh$, and $T_1 \stb T_k$ are right triangles
having perpendicular sides parallel to the axes.
\end{lemma}
\proof It is proved in pages 81-85 of Boltianskii's book \cite{B} that every
polygon $A$ is equidecomposable using only translations to the union of
nonoverlapping trapezoids having horizontal bases and a vertical leg.
One can easily check that the construction actually proves
$S$-equidecomposability. We sketch the argument.

The horizontal lines going through the vertices of $A$ dissect $A$ into
trapezoids and triangles such that each triangle has a horizontal side, and
the bases of the trapezoids are also horizontal. It is easy to see that the
triangles are $S$-equidecomposable to trapezoids with horizontal bases. (See Figure 49 on
page 82 of Boltianskii's book \cite{B}.)
Also, each trapezoid is $S$-equidecomposable to a union of right trapezoids
having horizontal bases and a vertical leg. The other legs of these trapezoids
are parallel to one of the sides of $A$.

Since each of the trapezoids obtained can be dissected into a rectangle
with sides parallel to the axes and a right triangle having perpendicular
sides parallel to the axes, the statement of the lemma
follows. \hfill $\square$ 

\begin{lemma}\label{l7}
Suppose $0,1\in S$.
Let $A,B\in \pp _S$, and suppose that $\nu _u (A)=\nu _u (B)$ for every unit
vector $u$. Then there are nonoverlapping polygons $H_1 , H_2 \in \hh$ and
$C \in \pp _S$ such that $A\sim _S H_1 \cup C$ and $B\sim _S H_2 \cup C$.
\end{lemma}
\proof
Let $U_S$ denote the set of unit vectors $u=(x,y)$ such that $x\ne 0$, $y\ne 0$
and $y/x\in S$. If $u\in U_S$, then we denote by $\ttt _{u}^+$ (resp. $\ttt _u^-$)
the set of right triangles $T$ with perpendicular sides parallel to the axes
and such that their hypotenuse is parallel to $u$, and $\nu _u (T)>0$ (resp.
$\nu _u (T)<0$).

By Lemma \ref{l7a}, $A\sim _S H\cup T_1 \cup \ldots \cup T_k$ and
$B\sim_S H' \cup T'_1 \cup \ldots \cup T'_n$, where $H, H'\in \hh$, and
$T_1 \stb T_k ,T'_1 \stb T'_n \in \pp _S$ are right triangles with
perpendicular sides parallel to the axes.

For every $u\in U_S$ we have, by assumption, $\nu _u (A)= \nu _u (B)$.
Let $a_u$ be the common value. 

For a given $u\in U_S$, we can translate the triangles $T_i$ belonging to
$\ttt_{u}^+$ such that the translated copies are nonoverlapping, and the
union of their hypotenuses is a segment $I_u$. Let $D_u$ denote the union of
these translated triangles. Note that $D_u$ lies either below or above the
segment $I_u$, depending on $u$. Similarly, we can
translate the triangles $T_i \in \ttt _u^-$ such that the union of their
hypotenuses is a segment $J_u$. Let $E_u$ denote the union of these
translated triangles. Clearly, $a_u =\nu _u (A)=|I_u |-|J_u |$.

If $a_u =0$ then, translating $E_u$, we may assume that $I_u =J_u$. Note that
$D_u$ and $E_u$ lie in different sides of $I_i$, and $H_u =D_u \cup E_u \in \hh$.
Since $\nu _u (B)=a_u =0$, a similar construction shows that suitable translated
copies of the triangles $T'_j$ belonging to $\ttt_{u}^+ \cup \ttt_{u}^-$ is a
polygon $H'_u \in \hh$.

If $a_u >0$ and $I_u =[a_u ,b_u ]$, then we may assume that $J_u =[a_u ,c_u ]$,
where $c_u \in I_u$. It is easy to check that $D_u \cup E_u = H_u \cup G_u$,
where $H_u \in \hh$, and $G_u$ is the union of some triangles 
such that the union of their hypotenuses is $[c_u ,b_u ]$.

Since $a_u =\nu _u (B)$, we can translate the triangles $T'_j$ belonging to
$\ttt_{u}^+ \cup \ttt_{u}^-$ such that their union  equals $H'_u \cup G'_u$,
where $H'_u \in \hh$, and $G'_u$ is the union of some triangles 
such that the union of their hypotenuses is $[c_u ,b_u ]$.

Note that $G_u$ and $G'_u$ lie on the same side of $[c_u ,b_u ]$.
Put $C_u =G_u \cap G'_u$. It is easy to check that $G_u \se C_u \in \hh$
and $G'_u \se C_u \in \hh$. 

If $a_u <0$, then we have a similar construction with the roles of $I_u$ and
$J_u$ exchanged.

One can easily see that the sets
$$H_1 =H\cup \bigcup_{u\in U_S} H_u \cup \bigcup_{u\in U_S , a_u \ne 0} (G_u \se C_u ),
\ H_2 =H'\cup \bigcup_{u\in U_S} H'_u \cup \bigcup_{u\in U_S , a_u \ne 0} (G'_u \se C_u )$$
and $C=\bigcup_{u\in U_S , a_u \ne 0} C_u$ satisfy the requirements. \hfill $\square$

\section{Proof of Theorem \ref{t4}} \label{s4}
Let $K$ denote the field generated by $S$. (If $S=\{ 0\}$, then $K=\qq$.)
Let $A,B\in \pp _S$ be given. If $A\sim _S B$, then $\mu (A)=\mu (B)$ holds
for every invariant on $\pp _S$. If $\mu$ is an invariant of $\pp _K$, then
the restriction of $\mu$ to $\pp _S$ is automatically an invariant on
$\pp _S$, and thus we have $\mu (A)=\mu (B)$.

Next suppose that $\mu (A)=\mu (B)$ whenever $\mu$ is an invariant of the first
or of the second kind on $\pp _S =\hh$ or on $\pp _K$ according to the cases
$S=\{ 0\}$ and $\{ 0,1\} \su S$. We have to show that $A\sim _S B$.

We have, by assumption, $\nu _u (A)=\nu _u (B)$ for every unit vector $u$.
By Lemma \ref{l7}, we may assume that $A=H_1 \cup C$ and $B=H_2 \cup C$,
where $H_1 ,H_2, C$ are nonoverlapping polygons, $H_1 , H_2 \in \hh$ and
$C \in \pp _S$. Clearly, it is enough to show that $H_1 \sim _S H_2$.
If $\mu$ is an invariant of the second kind and is defined on either $\hh$
or on $\pp _K$ according to the cases $S=\{ 0\}$ and $0,1\in S$, then
we have
$$\mu (H_1 )+\mu (C)=\mu (A)=\mu (B)=\mu (H_2 )+\mu (C),$$
and thus $\mu (H_1 )=\mu (H_2 )$.

We have $H_1 = \bigcup_{i=1}^p R_i$ and $H_2 =\bigcup_{j=1}^r Q_j$,
where the systems $\{ R_1 \stb R_p \}$, $\{ Q_1 \stb Q_r \}$ consist of
nonoverlapping rectangles with sides parallel to the axes. Let
the lengths of the sides of $R_i$ be $a_i$ and $b_i$ $(i=1\stb p)$, and those
of $Q_j$ be $c_j$ and $d_j$ $(j=1\stb r)$.
By Lemma \ref{l6}, there are positive numbers $h_1 \stb h_t$ such that they
are linearly independent over $K$, and each of the numbers $a_i ,b_i ,c_j ,d_j$
is a linear combination of $h_1 \stb h_t$ with nonnegative coefficients
belonging to $K$. Using suitable vertical and horizontal lines we can decompose
the rectangles $R_1 \stb R_p , Q_1 \stb Q_r$ into rectangles of size
$\al h_i \times \be h_j$, where $\al , \be$ are positive elements of $K$.
We have  
$$([0, \al h_i ]\times [0,\be h_j ])\sim _S ([0, h_i ]\times [0,\al \be h_j ]).$$
Indeed, if $S=\{ 0\}$, then this follows from (i) of Lemma \ref{l4}, and if
$0,1\in S$, then from Theorem \ref{t5}.
Thus $H_1$ is $S$-equidecomposable to the union of nonoverlapping rectangles
with sides parallel to the axes, and of size $h_i \times \ga h_j$, where
$\ga \in K$. If, among these rectangles, there are more than one
with the same pair $(i,j)$, then placing them on the top of each other,
we unify them into one single rectangle.

Summing up: there is a polygon $D_1$ such that $H_1 \sim _S D_1$, and
$D_1 =\bigcup_{(i,j)\in I} R_{i,j}$, where
$I\su \{ (i,j)\colon 1\le i, j \le t\}$ and $R_{i,j}$ is a rectangle
with sides parallel to the axes, and of size $h_i \times \ga_{i,j} h_j$, where
$\ga_{i,j}\in K$. Similarly, we find a polygon $D_2$
such that $H_2 \sim _S D_2$, and $D_2 =\bigcup_{(i,j)\in J} Q_{i,j}$, where
$J\su \{ (i,j)\colon 1\le i, j \le t\}$ and $Q_{i,j}$ is a rectangle
with sides parallel to the axes, and of size $h_i \times \de _{i,j} h_j$, where
$\de_{i,j}\in K$.
If $\mu$ is an arbitrary invariant on $\hh$ or on $\pp _K$, then we have
\begin{equation}\label{e8}
  \sum_{(i,j)\in I} \mu (R_{i,j})=\mu (D_1 )=\mu (H_1 )=\mu (H_2 )=\mu (D_2 )=
  \sum_{(i,j)\in J} \mu (Q_{i,j}).
\end{equation}
Now we consider the cases $S=\{ 0\}$ and $0,1\in S$ separately.
Suppose first $S=\{ 0\}$.
Let $x_1 \stb x_t , y_1 \stb y_t$ be arbitrary real numbers. Since
$h_1 \stb h_t$ are linearly independent over $\qq$, there are additive
functions $f,g\colon \rr \to \rr$ such that 
$f(h_i )=x_i$ and $g(h_i )=y_i$ $(i=1\stb t)$. 

Putting $\mu ([a,b]\times [c,d])=f(b-a)\cd g(d-c)$ we define an additive
function defined on the set of rectangles with sides parallel to the axes.
It is easy to check that $\mu$ can be extended to an invariant on $\hh$.
Thus $\mu$ is an invariant of second type on $\hh$ and then,
by \eqref{e8}, we obtain
$$\sum_{(i,j)\in I} \ga _{i,j} \cd x_i y_j =\sum_{(i,j)\in J} \de _{i,j} \cd x_i y_j .$$
Since $x_1 \stb x_t , y_1 \stb y_t$ were arbitrary, we find
that the polynomials \break
$\sum_{(i,j)\in I} \ga _{i,j} \cd x_i y_j$ and $\sum_{(i,j)\in J} \de _{i,j}
\cd x_i y_j$ are identical; that is, $I=J$ and $\ga _{i,j} = \de _{i,j}$
for every $(i,j)\in I$. Thus the rectangles $R_{i,j}$ can be translated into
the rectangles $Q_{i,j}$, proving $D_1 \sim _S D_2$. Therefore, we have
$A\sim _S H_1 \cup C_1 \sim _S H_2 \cup C_2 \sim _S B$. This completes the
proof in the case $S=\{ 0\}$.

Next suppose $0,1\in S$. If $(i,j)\in I$ and $i>j$, then, by (iii) of
Lemma \ref{l4} and by Theorem \ref{t5}, 
$$([0, h_i ]\times [0, \ga_{i,j} h_j ]) \sim _S ([0, \ga_{i,j} h_j ]\times
[0, h_i ])\sim _S  ([0, h_j ]\times [0, \ga_{i,j}   h_i ]).$$
Then we can replace the rectangles $R_{i,j}$ with $i>j$ by rectangles
$R'_{j,i}$. Therefore, we may assume that $I\su \{ (i,j)\colon 1\le i\le j
\le t\}$ and, similarly, $J\su \{ (i,j)\colon 1\le i\le j \le t\}$.

Let $x_1 \stb x_t$ be arbitrary real numbers. Since $h_1 \stb h_t$ are
linearly independent over $K$, there is a linear map $f\colon \rr \to \rr$
from the linear space of $\rr$ over the field $K$ into itself such that
$f(h_i )=x_i$ $(i=1\stb t)$. By Corollary \ref{c2}, $\mu _f$ is an invariant
on $\pp _K$ of the second type. Then \eqref{e8} gives
$$\sum_{(i,j)\in I} \ga _{i,j} \cd x_i x_j =\sum_{(i,j)\in J} \de _{i,j} \cd x_i x_j$$
for every $x_1 \stb x_t \in \rr$. Thus the polynomials 
$\sum_{(i,j)\in I} \ga _{i,j} \cd x_i x_j$ and $\sum_{(i,j)\in J} \de _{i,j}
\cd x_i x_j$ are identical. Since $I,J\su \{ (i,j)\colon 1\le i\le j \le t\}$,
this implies $I=J$ and $\ga _{i,j} = \de _{i,j}$
for every $(i,j)\in I$. Thus the rectangles $R_{i,j}$ can be translated into
the rectangles $Q_{i,j}$, proving $D_1 \sim _S D_2$. Therefore, we have
$A\sim _S H_1 \cup C \sim _S H_2 \cup C \sim _S B$. This completes the
proof. \hfill $\square$

\section{Squares}\label{s5}
In this section we prove Theorem \ref{t6}. Put $R=[0,a]\times [0,b]$.

(i) If $R \sim _S ([0,c]\times [0,c])$, then $c=\sqrt{ab}$.
By Theorem \ref{t5}, $([0,a]\times [0,b]) \sim _S ([0,\sqrt{ab}]\times
[0,\sqrt{ab}])$ holds if and only if $\sqrt{b/a}\in K$. \hfill $\square$

(ii) Suppose $R\sim _S Q$, where $Q$ is a square. We may assume that the sides
of $Q$ are not parallel to the axes, because otherwise the statement follows
from (i). Translating $Q$ we may also assume that the
points $(c,0)$ and $(0,d)$ are vertices of $Q$, where $c,d>0$. Then
$c^2 +d^2 =\la _2 (Q)=\la _2 (R)=ab$, and thus $a/b=\ga ^2 +\de ^2$, where
$\ga =c/b$ and $\de =d/b$.

We prove $\ga , \de \in K$. Suppose $\ga =c/b\notin K$. Then there is a basis
$B$ of $\rr$ as a linear space over $K$ such that $b,c \in B$. Let $f(x)$
denote the coefficient of $c$ in the representation of $x$ as a linear
combination of elements of $B$ with coefficients from $K$. Then
$f$ is additive, and $f(sx)=s\cd f(x)$ holds for every $x\in \rr$ and $s\in S$.
Then $\mu _f$ is an invariant on
$\pp _S$ by Corollary \ref{c2}, and thus $\mu _f (R)=\mu _f (Q)$.
It is easy to check that $\mu _f (R)=2\cd f(a)\cd f(b)$.

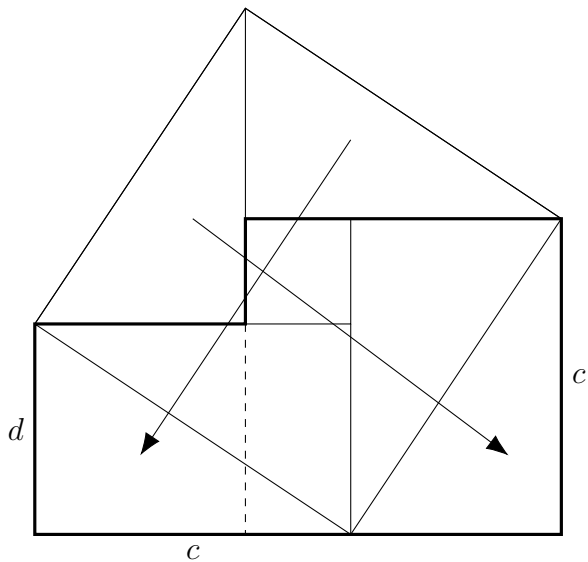
\begin{figure}[h!]
\centering
\begin{tikzpicture}[scale=0.7]

\coordinate (A) at (0,0);
\coordinate (B) at (10,0);
\coordinate (C) at (10,10);
\coordinate (D) at (0,10);

\coordinate (C_2) at (4,10);
\coordinate (C_1) at (10,6);

\coordinate (M_1) at (6,0);
\coordinate (M_2) at (0,4);
\coordinate (M_3) at (4,10);
\coordinate (M_4) at (10,6);

\coordinate (N_1) at (6,4);
\coordinate (N_2) at (4,4);
\coordinate (N_3) at (4,6);
\coordinate (N_4) at (6,6);

\coordinate (M) at (3,6);
\coordinate (N) at (9,1.5);
\coordinate (P) at (4,0);
\coordinate (Q) at (6,7.5);
\coordinate (R) at (2,1.5);

\draw (A) -- (B) -- (C_1) -- (C_2) -- (M_2) -- cycle;
\draw (M_1) -- (M_2) -- (M_3) -- (M_4) -- cycle;
\draw (N_1) -- (N_2) -- (N_3) -- (N_4) -- cycle;
\draw (N_1)--(M_1);
\draw (N_2)--(M_2);
\draw (N_3)--(M_3);
\draw (N_4)--(M_4);
\draw[very thick] (A) -- (M_2) -- (N_2) -- (N_3) -- (C_1) -- (B) -- cycle;

\draw[dashed] (N_2) -- (P);

\draw[-{Latex[length=3mm]}] (M) -- (N);
\draw[-{Latex[length=3mm]}] (Q) -- (R);


\node[left] at (0, 2) {$d$};
\node[below] at (3,0) {$c$};
\node[right] at (10, 3) {$c$};

\end{tikzpicture}

\caption{Equidecomposition of a square into two squares, $c=\gamma a$, $d=\delta b$}
    \label{fig2}
\end{figure}

Since $Q\in \pp _S$, the slopes of the sides of $Q$ belong to $S$. As Figure
\ref{fig2} shows\footnote{Figure \ref{fig2} comes from Airy's proof of the
Pythagoras theorem \cite[p. 4]{F}.}, $Q$ is $S$-equidecomposable to the union of two squares of side length
$c$ and $d$, and thus $\mu _f (Q)=2\cd f(c)^2 +2\cd f(d)^2$. Therefore,
$\mu _f (R)=\mu _f (Q)$ gives $f(a)\cd f(b)= f(c)^2 +f(d)^2$.
However, we have $f(c)=1$ and $f(b)=0$, which is a contradiction.
Thus $\ga \in K$, and a similar argument shows $\de \in K$. This proves
the `only if' part of (ii).

To prove the `if' part, let $a=\ga ^2 b +\de ^2 b$, where $\ga , \de \in K$.
We have $R=R_1 \cup R_2$, where $R_1 =[0,\ga ^2 b]\times [0,b]$ and $R_2 =
[\ga ^2 b,a]\times [0,b]$. Now $(\ga ^2 b)/(\ga b)
=\ga \in K$ implies, by Theorem \ref{t5}, that $R_1 \sim _S Q_1$, where
$Q_1 =[0,\ga b ]\times [0,\ga b ]$. Similarly, $R_2 \sim _S Q_2$, where
$Q_2 =([0,\de b ]\times [0,\de b ]) \sim _S ([\ga b,\ga b+\de b] \times
[0,\de b ] )=Q'_2$. As Figure \ref{fig2} shows (with $c=\ga b$ and $d=\de b$),
$Q_1 \cup Q'_2$ is $K$-equidecomposable to a square, and then so is $R$.
By (ii) of Theorem \ref{t4}, this implies that $R$ is $S$-equidecomposable to a
square. \hfill $\square$

\section{Description of the invariants of $\pp _S$}\label{s6}
The invariants $\nu _u$ and $\mu _F$ were introduced in \eqref{ev} and in
Section \ref{s1}. In this section our aim is to describe {\it all}
invariants on $\pp_S$. 

If $S=\{ 0\}$, then $\pp _S =\hh$. If $\mu$ is an invariant of $\hh$, then
putting
\begin{equation}\label{e9}
F(x,y)=\mu ([0,x] \times [0,y])
\end{equation}
we define a function $F$ mapping $\{ (x,y)\colon x,y >0\}$ into $\rr$.
Since $\mu$ is translation invariant and additive, we have
\begin{equation}\label{e9a}
\begin{split}
  F(x_1 +x_2 ,y)=& \mu ([0,x_1 +x_2 ] \times [0,y])=\\
&  \mu ([0,x_1 ] \times [0,y])
+\mu ([x_1 ,x_1 +x_2 ] \times [0,y])=\\
& F(x_1 ,y)+ F(x_2 ,y)
\end{split}
\end{equation}
for every $x_1 ,x_2 ,y>0$. Similarly, we have $F(x,y_1 +y_2 )=F(x,y_1 )+
F(x,y_2 )$ for every $x,y_1 ,y_2 >0$. This easily implies that $F$ can be
extended to $\sik$ as a biadditive function. Clearly, if $H\in \hh$ and
$H$ is the union of the nonoverlapping rectangles $[a_i ,b_i ]\times [c_i ,d_i ]$, then
\begin{equation}\label{e10}
  \mu (H)=\sumin F(b_i -a_i ,d_i -c_i ).
\end{equation}
In the other direction, if $F\colon \sik \to \rr$ is a biadditive function, then
\eqref{e10} defines an invariant on $\hh$. In this way we have described the
invariants of $\hh$.

Now let $0,1\in S$. Let $U$ denote the set of unit vectors. As we saw
before, $\nu _u$ is an invariant on the set $\pp$ of all polygons for every
$u \in U$. Also, if $f\colon \rr \to \rr$ is additive, then $f\circ \nu _u $
is an invariant on $\pp$ as well. Let $f_u \colon \rr \to \rr$ be an
additive function for every unit vector $u$, and put
\begin{equation}\label{e11}
\nu (A)=\sum _{u\in U} f_u (\nu _u (A))
\end{equation}
for every polygon $A$. The sum in the right hand makes sense, since all but a
finite number of terms vanish. It is clear that $\nu$ is an invariant of $\pp$
for every choice of the additive functions $f_u$ $(u\in U)$.
Note that the invariant $\nu$ defined in \eqref{e11} has the property that
$\nu (t\cd A)=t\cd \nu (A)$ for every polygon $A$ and for every positive
rational number $t$. That is, these invariants are rationally $1$-homogeneous.
By \cite[Theorem 4.1]{KP}, every rationally $1$-homogeneous
invariant on $\pp$ can be obtained this way.

\begin{theorem}\label{t8}
Suppose that $0,1\in S$, and let $\mu$ be an invariant on $\pp _S$.
Then $\mu$ is the restriction to $\pp _S$ of a rationally $1$-homogeneous
invariant on $\pp$ if and only if $\mu(H)=0$ for every $H\in \hh$.
\end{theorem}
\proof The `only if' direction is clear, so it is enough to prove the other
direction. Suppose $\mu$ is an invariant on $\pp _S$ vanishing on $\hh$.
We show that $\mu$ is the restriction to $\pp _S$ of a rationally
$1$-homogeneous invariant on $\pp$.

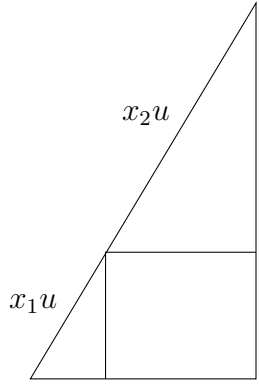
\begin{figure}[h!]
\centering
\begin{tikzpicture}[scale=0.5]

\coordinate (A) at (0,0);
\coordinate (B) at (6,0);
\coordinate (C) at (2,0);
\coordinate (D) at (6,3.37);
\coordinate (E) at (6,10);
\coordinate (F) at (2,3.37);

\draw (A) -- (B) -- (E) --  cycle;
\draw (C) -- (F) -- (D);

\node[left] at (1, 2) {$x_1u$};
\node[left] at (4,7) {$x_2u$};

\end{tikzpicture}

\caption{Additivity of $f_u$}
    \label{fig3}
\end{figure}

Let $U_S$ denote the set of unit vectors $u=(x,y)$ such that $x\ne 0$,
$y\ne 0$ and $y/x\in S$. Let $u\in U_S$ be fixed, and denote by $T^u_x$ the
right triangle with the following properties: its perpendicular sides are
parallel to the axes, and $\nu _u (T^u_x )=x$. (In particular, its hypotenuse is parallel to $u$ and has length $x$.) Putting $f_u (x)=\mu (T^u_x )$ for every
$x>0$, we define a function $f_u \colon (0,\infty )\to \rr$. If $x_1 ,x_2 >0$,
then $T^u_{x_1 +x_2}$ is the union of the triangles $T^u_{x_1}$ and $T^u_{x_2}$, and
a rectangle (see Figure \ref{fig3}). Since $\mu$ vanishes on $\hh$, we obtain
$f(x_1 +x_2 )=f(x_1 )+f(x_2 )$ for every $x_1 ,x_2 >0$. We can extend $f_u$
to $\rr$ as an additive function, also denoted by $f_u$. Now we put
$\nu (A)=\sum _{u\in U_S} f_u (\nu _u (A))$ for every $A\in \pp$. Clearly,
we have $\nu (T^u_x )=\mu (T^u_x )$ for every $u\in U_S$ and $x>0$.

We also have $\nu (H)=\mu (H)=0$ for every $H\in \hh$. Therefore, by Lemma
\ref{l7a}, we have $\nu (A)=\mu (A)=0$ for every
$A\in \pp _S$. \hfill $\square$
\begin{remark}\label{r3}
{\rm The proof above shows that the representation \eqref{e11} of the
rationally $1$-homogeneous invariant on $\pp _S$ is not unique. Let
$\nu _{(1,0)}$ be the invariant of first kind corresponding to the horizontal
unit vector $(1,0)$. The proof of Theorem \ref{t8} shows that $\nu _{(1,0)}$
can be represented in the form $\sum _{u\in U \se \{ (1,0) \} } f_u \circ \nu _u$.
}
\end{remark}

The following theorem gives the Hadwiger-McMullen decomposition of the
invariants on $\pp _S$.

\begin{theorem}\label{t9}
Suppose that $0,1\in S$. Then every invariant $\mu$ on $\pp _S$ can
be represented uniquely as a sum $\nu +\mu _F$, where $\nu$ is the restriction
to $\pp _S$ of a rationally $1$-homogeneous invariant on $\pp$, and $F\colon
\sik \to \rr$ is a symmetric, biadditive function satisfying \eqref{e1}.
\end{theorem}
\proof
Let $\mu$ be an invariant on $\pp _S$. Then \eqref{e9} defines a
function $F$ on $(0,\infty )\times (0,\infty )$, and can be extended to $\sik$
as a biadditive function (see \eqref{e9a}). By Lemma \ref{l4}, we have
$([0,x]\times [0,y]) \sim _S ([0,y]\times [0,x])$ for every $x,y>0$. Therefore,
we have $F(x,y)=F(y,x)$ if $x,y>0$. Since $F$ is biadditive on $\sik$,
we have $F(x,y)=F(y,x)$ for every $x,y$; that is, $F$ is symmetric.

If $s\in S$ and $s\ne 0$ then, by Theorem \ref{t5}, we have
$([0,|s|\cd x]\times [0,y]) \sim _S ([0,x]\times [0,|s|\cd y])$ for every
$x,y>0$.
Thus $F(|s| \cd x,y)=F(x,|s|\cd y)$ for every $x,y>0$, $s\in S$, $s\ne 0$.
Since $F$ is biadditive on $\sik$, $F(sx,y)=F(x,sy)$ for every $x,y$ and for
every $s\in S$; that is, $F$ satisfies \eqref{e1}. Thus $\mu _F$ is an
invariant on $\pp _S$ by Theorem \ref{t7}.

Let $\nu =\mu -\mu _F$. Then $\nu$ is an invariant on $\pp _S$, and vanishes on
$\hh$. By Theorem \ref{t8}, $\nu$ is the restriction to $\pp _S$ of a
rationally $1$-homogeneous invariant on $\pp$.

The uniqueness of the representation follows from the fact that if an
invariant is rationally $1$-homogeneous and rationally $2$-homogeneous
simultaneously, then it is zero. \hfill $\square$

\begin{cor}\label{c4}
If $\mu$ is a rationally $2$-homogeneous invariant on $\pp _S$, then
$\mu =\mu _F$ for a suitable symmetric, biadditive function $F\colon
\sik \to \rr$ satisfying \eqref{e1}.
\end{cor}

\proof
By Theorem \ref{t9}, $\mu =\nu +\mu _F$, where $\nu$ is rationally
$1$-homogeneous, and $f$ is a symmetric, biadditive function 
satisfying \eqref{e1}. For every $A\in \pp _S$ we have
$$t^2 \cd \mu (A)=\mu (tA)=\nu (tA)+\mu _F (tA)= t\nu (A)+t^2 \mu _F (A)$$
for every positive rational $t$. This implies $\nu (A)=0$ and $\mu (A)=
\mu _F (A)$. \hfill $\square$

\section{Duality}
We conclude with some remarks concerning the duality between the subsets of
$\rr$ and the sets of symmetric biadditive functions. As we saw above, the
rationally $1$-homogeneous invariants are defined and are invariants on the
set of all polygons.  On the other hand, if $0,1\in S$ and if $F$ is a
symmetric and biadditive function, then, although $\mu _F$ is also
defined on all polygons, {\it $\mu _F$ is an invariant on $\pp _S$ if and only
if $F$ satisfies \eqref{e1}.} As for the `only if' part: by Theorem \ref{t5},
the rectangles $[0,tx]\times [0,y]$ and $[0,x]\times [0,t y]$ are
$S$-equidecomposable for every $x,y,t>0$, where $t=|s|$ for some $s\in S$.
Therefore, if $\mu _F$ is an invariant on $\pp_S$, then
$$2F(tx,y)=\mu _F ([0,t x]\times [0,y])=\mu _F ([0,x]\times
[0,ty]) = 2F(x, t y).$$
Since $F$ is biadditive, this implies that $F(sx,y)=F(x,sy)$ 
for every $x,y\in \rr$ and $s\in S$.

Let $\ff$ denote the family of all symmetric and biadditive
functions on $\sik$. The observation above suggests that a certain
duality exists between subsets of $\rr$ and subsets of $\ff$.
If $F\in \ff$, then we put
$$F^\perp =\{ s\in \rr \colon F(sx,y)=F(x,sy) \ (x,y\in \rr )\} .$$
From the considerations above it follows that if $0,1\in S$, then 
{\it $\mu _F$ is an invariant on $\pp _S$ if and only if $S\su F^\perp$.}
\begin{proposition}\label{p4}
\begin{enumerate}[{\rm (i)}]
\item The set $F^\perp$ is a subfield of $\rr$ for every $F\in \ff$.
\item For every subfield $K$ of $\rr$ there is an $F\in \ff$ such that
$F^\perp =K$.
\end{enumerate}
\end{proposition}
\proof (i): It easily follows from the symmetry and the biadditivity of $F$ that
$\qq \su F^\perp$ and that $F^\perp$ is an additive subgroup of $\rr$.
If $s,t\in F^\perp$, then we have $F((st)x,y)=F(tx,sy)=F(x,(st)y)$ for every
$x,y$, and thus $st\in F^\perp$. Finally, if $s\in F^\perp$ and $s\ne 0$, then
$$F(x/s,y)=F(y,x/s)=F(s\cd (y/s), x/s)=F(y/s,x)=F(x,y/s)$$
for every $x,y$, and thus $1/s\in F^\perp$.

(ii) Let $K$ be a subfield of $\rr$; then $\rr$ is a linear space over $K$.
Let $B$ be a basis of this linear space such that $1\in B$, and let $g(x)$
denote the coefficient of $1$ in the representation of $x$ as a linear
combination of elements of $B$ with coefficients from $K$. Then 
$g$ is an additive function such that $g(sx)=s\cd g(x)$ for every $s\in K$
and $x\in \rr$. Also, we have $g(x)=x\iff x\in K$ for every $x$.

Put $f(x)=g(x)-x$ $(x\in \rr )$. Then $f$ is an additive function, $f(sx)
=s\cd f(x)$ for every $s\in K$ and $x\in \rr$, and $f(x)=0\iff x\in K$ for
every $x\in \rr$. Let $F(x,y)=f(x)\cd f(y)$ for every $x,y\in \rr$. It is clear
that $F\in \ff$, and $K\su F^\perp$.

Suppose $t\in F^\perp$; we prove $t\in K$. Since $t\in F^\perp$, we have
$$f(t\cd 1)\cd f(t)=F(t\cd 1,t)=F(1,t^2 )=f(1)\cd f(t^2 )=0.$$
Thus $f(t)^2 =0$, $f(t)=0$ and $t\in K$. \hfill $\square$

\begin{remark}
{\rm 
By (i) of the Proposition \ref{p4} we can see that {\it the maximal $S$ such
that $\mu _F$ is an invariant on $\pp _S$ is always a field.}
}
\end{remark}

\begin{cor}\label{c3}
Suppose $0,1\in S$, and let $K$ denote the field generated by $S$. If $S'
\not\subsetneq K$, then there are rectangles $A,B$ with sides
parallel to the axes such that $A\sim _{S'} B$, but $A\sim_S B$ does not hold.
\end{cor}
\proof Let $t\in S' \se K$. By (ii) of Proposition \ref{p4}, there is an $F\in \ff$ such that
$F^\perp =K$. Since $t\notin K$, this implies that $F(tx,y)\ne F(x,ty)$ for some $x,y\in \rr$.
We may assume that $t,x,y>0$. Let $A=[0,tx]\times [0,y]$ and $B=[0,x]\times [0,ty]$. Then we have 
$A\sim _{S'} B$, since $t\in S'$ (see (ii) of Lemma \ref{l4}). On the other hand,
$A\sim_S B$ does not hold, since $\mu _F$ is an invariant on $\pp _S$ by Theorem \ref{t7},
and $\mu _F (A) =2F(tx,y)\ne 2F(x,ty)=\mu _F (B)$. \hfill $\square$

If $\gg \su \ff$, then we put $\gg ^\perp =\bigcap _{F\in \gg} F^\perp$.
Clearly, $\gg ^\perp$ is a field, and it is the maximal subset $S$ of $\rr$ such
that $\mu _F$ is an invariant on $\pp _S$ for every $F\in \gg$.

In the other direction, let
$$S^\perp =\{ F\in \ff \colon \mu _F \ \text{is an invariant on} \ \pp _S \} .$$
Or, equivalently, let $S^\perp =\{ F\in \ff \colon S\su F^\perp \}$.
\begin{proposition}\label{p5}
If $0,1\in S$, then $\left( S^\perp \right) ^\perp$ equals
the field generated by $S$. Thus $\left( S^\perp \right) ^\perp =S$
for every subfield $S$ of $\rr$.
\end{proposition}
\proof It is enough to prove the first statement. Let $K$ denote the field
generated by $S$. It is clear that $K\su \left( S^\perp \right) ^\perp$. By
Proposition \ref{p4}, there is an $F\in \ff$ such that $F^\perp =K$. Then
$F\in S^\perp$, and $\left( S^\perp \right) ^\perp \su F^\perp =K$. \hfill $\square$

Note that $S^\perp$ is always a linear subspace of $\ff$. Moreover,
$S^\perp$ has the following property: if a function $G\colon \sik \to \rr$
is such that for every finite set $X\su \sik$ there is an $F\in S^\perp$
with $G|_X =F|_X$, then $G\in S^\perp$. (This implies that $S^\perp$ 
is a closed subspace of the product space $\prod_{i\in \rr} Y_i$, where 
each $Y_i$ equals $\rr$ equipped with the discrete topology.)

We may ask whether $\left( \gg ^\perp \right) ^\perp =\gg$ holds at least in
those cases, when $\gg$ is a closed subspace of $\ff$. We show that the
answer is negative.

By (ii) of Proposition \ref{p4}, there is an $F\in \ff$ such that $F^\perp =\qq$.
Let $\gg =\{ c\cd F\colon c\in \rr \}$. It is easy to check that 
$\gg$ is a closed linear subspace of $\ff$. Now we have $\gg ^\perp =\qq$ and
$$\left( \gg ^\perp \right) ^\perp =\qq ^\perp =\ff \ne \gg .$$

\subsection*{Acknowledgments}
We thank Péter Medvegyev for bringing our attention to the problem of equidecompositions using dissections restricted to given directions.

Both authors were supported by the Hungarian National Foundation for
Scientific Research, Grant No. K146922. The first author
was also supported by the J\'anos Bolyai Research Fellowship and Hungarian National Foundation for
Scientific Research, Grant No. Starting 150576.

\begin{small}\noindent
(G. Kiss)\\
{\sc Corvinus University of Budapest, Department of Mathematics \\
Fővám tér 13-15, Budapest 1093, Hungary,\\
and\\
HUN-REN Alfred Renyi Mathematical Institute\\
Reáltanoda street 13-15, H-1053, Budapest, Hungary\\
E-mail: {\tt kigergo57@gmail.com}}

(M. Laczkovich)\\
{\sc E\"otv\"os Lor\'and University\\
Budapest, P\'azm\'any P\'eter s\'et\'any 1/C, 1117 Hungary\\
E-mail: {\tt miklos.laczkovich@gmail.com}}
\end{small}

\end{document}